\documentstyle[12pt,leqno,amssymb,graphics]{amsart}
\newtheorem{thm}{Theorem}[section]

\newtheorem{prop}[thm]{Proposition}
\newtheorem{cor}[thm]{Corollary}

\theoremstyle{definition}
\newtheorem{dfn}[thm]{Definition}

\theoremstyle{remark}

\begin{document}

\newcommand{\su}{\subseteq}
\newcommand{\pa}{\partial}

\newcommand{\x}{\times}
\newcommand{\I}{[0,1]}

\newcommand{\bbox}{\hfill}

\newcommand{\R}{{\Bbb R}}
\newcommand{\Z}{{\Bbb Z}}
\newcommand{\E}{{{\Bbb R}^3}}
\newcommand{\C}{{{\Bbb Z}/2}}

\newcommand{\Q}{F\times [0,1]}
\newcommand{\s}{\Sigma F}

\newcommand{\p}{{\pi_{1}}}
\newcommand{\pp}{{\pi_{2}}}
\newcommand{\ppp}{{\pi_{3}}}

\newcommand{\ep}{{\epsilon}}
\newcommand{\0}{\times \{ 0 \}}
\newcommand{\1}{\times \{ 1 \}}
\newcommand{\h}{{1\over 2}}
\newcommand{\f}{{\hat{F}}}
\newcommand{\se}{{7\over 8}}
\newcommand{\wi}{{\widetilde{{\Z}^n}}}
\newcommand{\im}{{Imm(F,\E)}}

\newcounter{numb}

\title{Finite Order $q$-Invariants of immersions of surfaces into 3-Space}
\author{Tahl Nowik}
\address{Department of Mathematics, Columbia University, New York, 
NY 10027, USA.} 
\email{tahl@@math.columbia.edu}
\date{March 21, 1999}

\begin{abstract}
Given a surface $F$, we are interested in
$\C$ valued invariants 
of immersions of $F$ into $\E$, which
are constant on each connected component
of the complement of the 
quadruple point discriminant in $\im$.
Such invariants will be called ``$q$-invariants.''
Given a regular homotopy class $A\su \im$, we 
denote by $V_n(A)$ the space of all $q$-invariants on $A$ of 
order $\leq n$. We show that if $F$ is orientable, then
for each regular homotopy class $A$ 
and each $n$,
$\dim ( \  V_n (A) / V_{n-1}(A) \  ) \leq 1$. 

\end{abstract}

\maketitle

\section{Introduction}
Let $F$ be a closed surface.
Let $\im$ denote the space of all immersions of $F$ into $\E$
and let $I_0\su\im$ denote the space of all generic immersions.

\begin{dfn}\label{df1}
A function $f:I_0 \to\C$ will be called a ``$q$-invariant'' 
if whenever
$H_t:F\to\E$ ($0\leq t \leq 1$)
is a generic regular homotopy with no quadruple points, then 
$f(H_0)=f(H_1)$. 
\end{dfn}

\begin{dfn}\label{df2}
Let $I_n \su \im$ denote the space of all immersions 
whose unstable self intersection consists
of precisely $n$ generic quadruple points,
and let $I=\bigcup_{n=0}^\infty I_n$. 
\end{dfn}

\begin{dfn}\label{df25}
Given a $q$-invariant $f:I_0\to\C$ we extend it to $I$ as follows:
For $i\in I_n$ let $i_1,...,i_{2^n}\in I_0$ be the $2^n$ generic immersions
that may be obtained by slightly deforming $i$.  Define
$$f(i)=\sum_{k=1}^{2^n} f(i_k).$$ 
\end{dfn}

For any $q$-invariant, we will always assume without mention
that it is extended to the whole of $I$ 
as in Definition \protect\ref{df25}.

The following relation clearly holds:

\begin{prop}\label{pr1}
Let $f$ be a $q$-invariant.
Let $i\in I_n$, $n\geq 1$, and let $p\in\E$ be one of its $n$ quadruple points.
Then: $f(i)=f(i_1)+f(i_2)$ where $i_1 , i_2 \in I_{n-1}$ are the two
immersions that may be obtained by slightly 
deforming $i$ in a small neighborhood of $p$.

(Or equivalently, since we are in $\C$,  $f(i_2)=f(i_1)+f(i)$.)
\end{prop}

\begin{dfn}\label{df3}
A $q$-invariant $f$ will be called 
``of finite order'' if $f|_{I_n}\equiv 0$ for some $n$. 

The ``order'' of
a finite order $q$-invariant $f$ 
is defined as the minimal $n$ such that 
$f|_{I_{n+1}}\equiv 0$.
\end{dfn}

(Compare our Definitions \protect\ref{df25} and \protect\ref{df3}
with 2.2 of [O].)

An example of a $q$-invariant of order 1 is the invariant
$Q$ which is defined by the property that if 
$H_t:F\to \E$ ($0\leq t \leq 1$) is a generic regular
homotopy in which $m$ quadruple points occur, 
then $Q(H_1)=Q(H_0)+m \mod 2$. 
In other words $Q$ is defined by the property that $Q|_{I_1} \equiv 1$.
It was proved in [N] 
that $Q$ indeed exists for any surface $F$.

There are $M=2^{2-\chi(F)}$ regular homotopy classes (i.e. connected
components) in $\im$.
Given a regular homotopy class $A\su\im$, we may repeat all our 
definitions with $A$ in place of $\im$. 
Let then $V_n(A)$ (respectively $V_n$)
denote the space of all $q$-invariants 
on $A$ (respectively $\im$) of order $\leq n$. 
$V_n(A)$ and $V_n$ are vector spaces over $\C$, and
$V_n = \bigoplus_{\alpha =1}^M V_n(A_\alpha)$ where
$A_1,...,A_M$ are the regular homotopy classes in $\im$. 
More precisely, a function $f:I_0\to\C$ is
a $q$-invariant of order $\leq n$  iff 
for every $1 \leq \alpha \leq M$,
$f|_{I_0 \cap A_\alpha}$ is a $q$-invariant of order $\leq n$.
And so studying $q$-invariants on $\im$ is the same as studying 
$q$-invariants on the various regular homotopy classes.

The purpose of this work is to prove the following:

\begin{thm}\label{th0}
If $F$ is orientable then 
$\dim ( \  V_n (A) / V_{n-1}(A) \   ) \leq 1$ 
for any $A$ and $n$.  
\end{thm}

By [N] $\dim ( \  V_1 (A) / V_0 (A)  \  ) \geq 1$ for any $A$
(for all surfaces, not necessarily orientable) and so we get:

\begin{cor}\label{cr1}

If $F$ is orientable then 
$\dim ( \  V_1 (A) / V_0 (A)  \  ) = 1$ for any $A$.

\end{cor}

Since as mentioned, $V_n = \bigoplus_{\alpha=1}^M V_n(A_\alpha)$, we get: 

\begin{cor}\label{cr2}
If $F$ is orientable of genus $g$ then 
$\dim (   V_n / V_{n-1}   ) \leq 2^{2g}$ 
for every $n$, and 
$\dim (   V_1 / V_0   ) = 2^{2g}$.

\end{cor}

\section{General $q$-Invariants}

The results in this section will not assume that the $q$-invariant
$f$ is of finite order.

\begin{thm}[The 10 Term Relation]\label{th1}
Let $i:F\to\E$ be any immersion whose non-stable self intersection
consists of one generic quintuple point, and some finite number of 
generic quadruple points. Let the quintuple point be located at
$p\in\E$ and let $S_1,...,S_5$ be the five sheets passing through $p$.
Let $i^1_k$ and $i^2_k$ ($k=1,...,5$) be the two immersions obtained from
$i$ by slightly pushing $S_k$ away from $p$ to either side. 
Then for any $q$-invariant $f$: 
$$\sum_{k=1}^5 \sum_{l=1}^2 f(i^l_k) = 0.$$
\end{thm}

\begin{pf}

Starting with $i$,
take $S_1$ and push it slightly to one side. Then take $S_2$ and 
push it away on a much smaller scale. What we now have is an immersion
$j$ where sheets $S_2,...,S_5$
create a little tetrahedron, and $S_1$ passes outside this tetrahedron.
We define the following regular homotopy $H_t:F\to \E$ 
beginning and ending with $j$,
we describe it in four steps: 
(a) $S_1$ sweeps to the other side of the tetrahedron. 
In this step four quadruple points occur.
(b) $S_2$ sweeps across the triple point of sheets $S_3,S_4,S_5$. This results
in the vanishing of the tetrahedron and its inside-out reappearance. 
One quadruple point occurs here. 
(c) $S_1$ sweeps back to its place. Four more quadruple points occur.
(d) $S_2$ sweeps back to its place. One more quadruple point occurs.

All together we have ten quadruple points, and say 
the $m$th quadruple point occurs at time $t_m$. 
It is easy to 
verify that 
the ten immersions
$H_{t_1},...,H_{t_{10}}$ are 
precisely (equivalent to) the ten immersions $i^l_k$ ($l=1,2$ , $k=1,...,5$.)
Also, $f(H_{t_m})=f(H_{{t_m}-\ep}) + f(H_{{t_m}+\ep})$
and so: $$\sum_{kl} f(i^l_k) = 
\sum_{m=1}^{10} f(H_{t_m})=
\sum_{m=1}^{10} (f(H_{{t_m}-\ep}) + f(H_{{t_m}+\ep})).$$ 
But $f(H_{{t_m}+\ep}) = f(H_{t_{m+1}-\ep})$ 
(where $m+1$ means $(m+1) \mod 10$)
and so this sum is 0.

\end{pf}

\begin{prop}\label{pr2}
Let $B(1)\su\E$ 
be the unit ball. Let $D_1(1),...,D_4(1) \su F$ be four disjoint discs 
which will each be parameterized as the unit disc, and let
$D(1)=\bigcup_{k=1}^4 D_k(1)$. 
Let $i\in I$ and
assume $i^{-1}(B(1))=D(1)$ and $i|_{D(1)}$
maps each $D_k(1)$ linearly onto some $L_k\cap B(1)$ where $L_k$ is a plane 
through the origin, and $L_1,...,L_4$ are in general position.
Let $i':D(1)\to B(1)$ be an immersion of the same sort as $i|_{D(1)}$
but with planes ${L'}_1,...,{L'}_4$.

For $0\leq r \leq 1$
let $B(r)\su B(1)$ and $D_k(r)\su D_k(1)$ 
be the ball and discs of radius $r$ and let $D(r)=\bigcup_{k=1}^4 D_k(r)$.

Then: There exists an immersion $j:F\to \E$ satisfying:
\begin{enumerate}
\item $j$ is regularly homotopic to $i$ via a regular homotopy
that fixes $F-D(1)$.
\item $j^{-1}(B(\h))=D(\h)$
\item $j|_{D(\h)}=i'|_{D(\h)}$
\item $f(j)=f(i)$ for any $q$-invariant $f$.
\end{enumerate}
\end{prop}

\begin{pf}
Slightly perturb $i$ if necessary so that the eight planes 
$L_k, {L'}_k$ will all be in general position. 
We define a regular homotopy $H_t$
from $i$ to an immersion $\tilde{i}$ as follows: 
Say $a$ is the point in $D_1(1)$ which is mapped to the origin.
Keeping $a$ and $F-D_1(1)$ fixed, we isotope $D_1(1)$ 
within $B(1)$ to get 
$\tilde{i}$ with ${\tilde{i}}^{-1}(B(\se))=D(\se)$ and 
$\tilde{i}|_{D_1(\se)}=i'|_{D_1(\se)}$. 

Let $i^1, i^2$ be the two
immersions obtained from $i$ by slightly pushing $D_1(1)$ off 
of the origin, and let
$\tilde{i}^1, \tilde{i}^2$ be the corresponding slight deformations
of $\tilde{i}$. 
$H_t$ induces regular homotopies $H^l_t$ ($l=1,2$)
from $i^l$ to $\tilde{i}^l$, 
and such that $H^l_t|_{D_1(1)}$ avoids the origin. 

Now, the only triple point of $\{ L_2, L_3, L_4 \}$ is the origin, 
and $H^l_t|_{D_1(1)}$ is an isotopy which avoids the origin,
and so $H^l_t$ will have no quadruple point,
and so $f(i^l)=f(\tilde{i}^l)$ ($l=1,2$). And so
(By Proposition \protect\ref{pr1})
$f(i)= f(i^1)+f(i^2) = f(\tilde{i}^1) + f(\tilde{i}^2) =f(\tilde{i})$.

We now repeat this process in the ball $B(\se)$ and with $D_2(\se)$,
obtaining an immersion $\tilde{\tilde{i}}$ with 
$\tilde{\tilde{i}}|_{D_1({6\over 8}) \cup D_2({6\over 8})} = 
i'|_{D_1({6\over 8}) \cup D_2({6\over 8})}$. After four iterations we 
get the desired $j$.

\end{pf}

\section{$q$-Invariants of Order $n$}

We now prove the following theorem, which clearly implies 
Theorem \protect\ref{th0} (our main theorem):

\begin{thm}\label{th2}
Assume $F$ is orientable and let $f$ be a
$q$-invariant of order $n$. 

Then for any regular 
homotopy class $A\su\im$, $f$ is constant on $I_n \cap A$.
\end{thm}

\begin{pf}

Let $i\in I$ and $p\in\E$ a quadruple point of $i$. A ball $B\su\E$
centered at $p$ as in Proposition \protect\ref{pr2},
i.e. such that $i^{-1}(B)$ is a union of four disjoint discs intersecting
in $B$ as four planes, will be called ``a good neighborhood for $i$ at $p$.''

For $i\in I_n$ 
let $p_1,...,p_n \in \E$ 
be the $n$ quadruple points of $i$ in some order,
and let $B_1,...,B_n$ be disjoint good
neighborhoods for $i$ at $p_1,...,p_n$.
We define $\pi_k(i) : F \to \pa B_k$ as follows:
Push each one of the four discs 
in $B_k$ slightly away from $p_k$ into the preferred side determined by 
the orientation of $F$. We now have a map that avoids $p_k$. 
Define $\pi_k(i)$ as the composition of this map with
the radial projection $\E -\{ p_k \} \to \pa B_k$. 

Let $d_k(i)$ denote
the degree of the map $\pi_k(i)$.

Let the symmetric group $S_n$ act on ${\Z}^n$ by
$\sigma (a_1,...,a_n)=(a_{\sigma(1)},...,a_{\sigma(n)})$, and let
$\widetilde{{\Z}^n}={\Z^n}/S_n$. Let the class
of $(a_1,...,a_n)$ in $\widetilde{{\Z}^n}$ be denoted by $[a_1,...,a_n]$.
For $i\in I_n$ we define $d(i)\in \widetilde{{\Z}^n}$ by
$d(i)=[d_1(i),...,d_n(i)]$. 

We break our proof into two steps. \emph{Step 1:}
If $i,j\in I_n \cap A$ and $d(i)=d(j)$ then $f(i) = f(j)$.
\emph{Step 2:} For any 
$(a_1,...,a_n)\in {\Z}^n$,
there are immersions $i,j \in I_n\cap A$ with
$d(i)=[a_1,a_2,...,a_n]$, $d(j)=[a_1 +1 , a_2,...,a_n]$
and $f(i)=f(j)$. The theorem clearly follows from these two claims.

\emph{Proof of Step 1:}
By composing $i$ with an isotopy $U_t : \E \to \E$ we may assume
that $p_1,...,p_n \in \E$ are the quadruple points of both $i$ and $j$
and that $d_k(i)=d_k(j)$ for each $1\leq k \leq n$. 
Let $B_1,...,B_n$ be disjoint good neighborhoods for 
both $i$ and $j$ at $p_1,...,p_n$.
By composing $i$ with an isotopy $V_t : F \to F$ we may further assume
that $i^{-1}(B_k) = j^{-1}(B_k)$ for every $k$.
We name the four discs in $F$ corresponding to $p_k$
by $D^{kl}$, $l=1,...,4$. 

Using Proposition \protect\ref{pr2} we may now change 
$i$ such that (for smaller $B_k$'s) we will have $i|_{D^{kl}}=j|_{D^{kl}}$
for all $1\leq k \leq n$ , $1 \leq l \leq 4$. The process of 
Proposition \protect\ref{pr2} indeed does not change $d_k(i)$, since
the slightly pushed discs appearing in the definition of $\pi_k(i)$
can follow the regular homotopy of Proposition \protect\ref{pr2} and this 
will induce a homotopy between the corresponding $\pi_k(i)$'s.

So we may assume $i|_{D^{kl}}=j|_{D^{kl}}$ 
for all $1\leq k \leq n$ , $1 \leq l \leq 4$.
We will now show that there exists a regular homotopy
from $i$ to $j$ such that each $D^{kl}$ moves only within its image in $\E$, 
and $F-\bigcup_{kl} D^{kl}$ moves only within 
$\E - \bigcup_k B_k$. We will then be done since such a regular homotopy
cannot change $f(i)$. Indeed,
no sheet will pass $p_1,...,p_n$ and so the only singularities 
that might be relevant are the quadruple points 
occurring in $\E-\bigcup_k B_k$.
But whenever such a quadruple point occurs, then we will have $n+1$
quadruple points all together, and so 
since $f$ is of order $n$, $f(i)$ will not change.
(Proposition \protect\ref{pr1}.)

To show the existence of the above regular homotopy,
we construct the following handle decomposition of $F$. 
Our discs 
$D^{kl}$, ($1\leq k \leq n$ , $1 \leq l \leq 4$) will be the 0-handles.
If $g$ is the genus of $F$ we will have $\ 2g+4n-1 \ $ 1-handles as follows:
$\ 2g \ $ 1-handles will have both ends glued to $D^{11}$ such that
$D^{11}$ with these $2g$ handles will decompose $F$
in the standard way. Then choose an ordering of
the discs $D^{kl}$ with $D^{11}$ first, and connect each two consecutive 
discs with a 1-handle. The complement of the 0- and 1-handles is one disc
which will be the unique 2-handle.

We first define our regular homotopy on the union of 0- and 1- handles.
Take a 1-handle $h$ of the first type. Since $i$ and $j$ 
are regularly homotopic, their restrictions to the annulus 
$D^{11} \cup h$ are also regularly homotopic. We can construct 
such a regular homotopy of $D^{11} \cup h$ fixing $D^{11}$ and avoiding
$\bigcup_k B_k$. 

Next consider the 1-handles of the second type. 
Take the 1-handle $h$ connecting
$D^{11}$ to the second disc in our ordering, call it $D'$. Then if
$i|_h$ and $j|_h$ are not regularly homotopic relative the gluing of $h$ to
$D^{11} \cup D'$, then we perform one full rotation of $D'$, 
as to make them regularly homotopic. (This will require a motion of the next
1-handle too.)  Again we perform all regular homotopies while
avoiding $\bigcup_k B_k$. We can now go along the chain of
1-handles of the second type, and regularly homotope them one by one as 
we did the first one. At each step we might need to move the next 0-handle and
1-handle, but we never need to change what we have already done.

So far we have constructed the desired regular homotopy on
the union of 0- and 1-handles. By [S] 
this regular homotopy may be extended to
the whole of $F$ (still avoiding $\bigcup_k B_k$.) 
And so, if we denote our 2-handle by $D$,
we are left with regularly homotoping $i|_D$ to $j|_D$ 
(relative $\pa D$.) Since $d_k(i)=d_k(j)$ for all $k$, these maps are
homotopic in $\E-\bigcup_k B_k$. 
By [S] they are also regularly homotopic in $\E-\bigcup_k B_k$, since the 
obstruction to that would lie in $\pi_2 (SO_3)=0$.

\emph{Proof of Step 2:} Take any immersion $i'\in I_n \cap A$ with
$d(i')=[a_1,...,a_n]$ and let $p_1,...,p_n \in \E$ be the 
quadruple points of $i'$, ordered such that $d_k(i')=a_k$, $1\leq k \leq n$. 
(Clearly any $[a_1,...,a_n]\in\wi$ may be realized 
within any regular homotopy class.) Take a disc in $F$ which is away from
the $p_k$'s and start pushing it (i.e. regularly homotoping it)
into its preferred side directing it towards
$p_1$. Avoid any of the $p_k$'s on the way, and so  
the immersion $i$ we will get just before arriving at $p_1$, will still
have $d_k(i)=a_k$ for all $k$.  We then pass $p_1$ creating a quintuple 
point, and continue to the other side arriving at an immersion $j$
which is again in $I_n$. Clearly $d_1(j)=a_1+1$ and 
$d_k(j)=a_k$ for $k\geq 2$.
We will now use Step 1 and the 10 term relation 
(Theorem \protect\ref{th1}) to show that $f(i)=f(j)$.
Indeed, let us name the five sheets of our quintuple point by 
$S_1,...,S_5$ where $S_1$ is the sheet coming from the disc that we pushed
into $p_1$. Let $i_m^1$ ($m=1,...,5$) denote the immersion obtained 
by pushing $S_m$ into its non-preferred side, and $i_m^2$ the immersion
obtained by pushing $S_m$ 
into its preferred side. Then $i=i_1^1$ and $j=i_1^2$.
Recall that $\pi_1(i_m^l)$ is constructed by pushing all four sheets involved 
in the quadruple point at $p_1$ into their preferred side. And so for each
$1\leq m\leq 5$, $\pi_1(i_m^1)$ has one sheet pushed into the non-preferred 
side and four sheets into the preferred side, and so $d_1(i_m^1)$ are all 
equal to each other. And, for each $1\leq m \leq 5$, $\pi_1(i_m^2)$ 
has all five sheets pushed into the preferred side and so
also $d_1(i_m^2)$ are all equal to each other. Clearly all this has no
effect on $d_k$ for $k\geq 2$, and so we have 
$d(i_m^1)=d(i)$ and $d(i_m^2)=d(j)$ for all $1\leq m\leq 5$.
And so by step 1, $f(i_m^1)=f(i)$ and $f(i_m^2)=f(j)$ for all $1\leq m\leq 5$.
And so by the 10 term relation,
$0=\sum_{ml} f(i_m^l) = 5 f(i) + 5 f(j) = f(i) + f(j)$ i.e. $f(i)=f(j)$.

\end{pf}

\end{document}